\documentclass[12pt]{amsart}
\usepackage{latexsym}
\usepackage{amssymb, amsmath,mathrsfs}
\usepackage{geometry}
\usepackage[OT2,T1]{fontenc}
\usepackage{setspace}
\usepackage{mathtools}

\usepackage[pagebackref,hypertexnames=false, colorlinks, citecolor=red, linkcolor=red]{hyperref}

\newtheorem{theorem}{Theorem}[section]
\newtheorem{lemma}[theorem]{Lemma}

\DeclareSymbolFont{cyrletters}{OT2}{wncyr}{m}{n}
\DeclareMathSymbol{\Sha}{\mathalpha}{cyrletters}{"58}

\theoremstyle{remark}
\newtheorem{remark}[theorem]{Remark}
\newtheorem{example}[theorem]{Example}

\begin{document}
\raggedbottom
\title{Characterizations of $A_2$ Matrix Power Weights}
\date{\today}

\author[K. Bickel]{Kelly Bickel$^{\dagger}$}
\address{Kelly Bickel, Department of Mathematics\\
Bucknell University\\
701 Moore Ave\
Lewisburg, PA 17837}
\email{kelly.bickel@bucknell.edu}
\thanks{$\dagger$ Research supported in part by National Science Foundation grant
DMS \#1448846.}

\author[K. Lunceford]{Katherine Lunceford}
\address{Katherine Lunceford, Department of Mathematics\\
Bucknell University\\
701 Moore Ave\
Lewisburg, PA 17837}
\email{}

\author[N. Mukhtar]{Naba Mukhtar}
\address{Naba Mukhtar, Department of Mathematics\\
Bucknell University\\
701 Moore Ave\
Lewisburg, PA 17837}
\email{}


\maketitle

\begin{abstract} In the scalar setting, the power functions $|x|^{\gamma}$, for $-1 < \gamma<1$, are the canonical examples of $A_2$ weights. In this paper, we study two types of power functions in the matrix setting, with the goal of obtaining canonical examples of $A_2$ matrix weights. We first study Type 1 matrix power functions, which are $n\times n$ matrix functions whose entries are of the form $a|x|^{\gamma}.$ Our main result characterizes when these power functions are $A_2$ matrix weights and has two extensions to Type $1$ power functions of several variables. We also study Type 2 matrix power functions, which are $n\times n$ matrix functions whose eigenvalues are of the form $a|x|^{\gamma}.$ We find necessary conditions for these to be $A_2$ matrix weights and give an example showing that even nice functions of this form can fail to be $A_2$ matrix weights.

\end{abstract}

\bibliographystyle{plain}

\section{Introduction}

\subsection{Background.}
Let $W: \mathbb{R} \rightarrow M_n(\mathbb{C})$ be an $n \times n$ matrix-valued function.
Then $W$ is a \emph{matrix weight} if $W(x)$ is positive definite for a.e.~$x \in \mathbb{R}$ and if each entry of $W$ is a locally integrable function. Such matrix weights are natural objects arising in a variety of settings, for example from spectral measures of stationary processes, from matrix Toeplitz operators, and in the multivariate case, in the study of systems of certain elliptic PDEs \cite{im16, tv97}. 
Given a matrix weight $W$, one can naturally consider the weighted space $L^2(W),$ which is the set of vector-valued functions $f: \mathbb{R} \rightarrow \mathbb{C}^n$ such that
\[ 
\| f\|^2_{L^2(W)} \equiv \int_{\mathbb{R}} \left \langle W(x) f(x), f(x) \right \rangle_{\mathbb{C}^n} dx =  \int_{\mathbb{R}} \left \| W(x)^{\frac{1}{2}} f(x) \right \|^2_{\mathbb{C}^n} dx <\infty.
\]
A standard question to ask about such spaces is: If an operator $T$ is bounded on $L^2(\mathbb{R}, \mathbb{C}^n)$, when is $T$ also bounded on $L^2(W)$? In \cite{tv97}, Treil and Volberg answered this question for the Hilbert transform and in \cite{nt97, vol97, cg01}, Nazarov and Treil, Volberg, and Christ and Goldberg separately answered the question for many important operators, including classes of Calder\'on-Zgymund operators and maximal functions. In all cases, the required condition is that $W$ be an $A_2$ matrix weight, namely
\begin{equation} \label{eqn:a2con} \big [ W \big ] _{A_2} \equiv \sup_{I} \left \| \langle W \rangle_I^{\frac{1}{2}}
 \langle W^{-1} \rangle_I^{\frac{1}{2}} \right \|^2 \approx \sup_{I} \text{Trace} \left( \langle W \rangle_I
 \langle W^{-1} \rangle_I \right)  < \infty,
\end{equation}
where the supremum is taken over all intervals $I$, $\langle W \rangle_I$ denotes the average $\frac{1}{|I|} \int_I W(x) dx$, and  $\left\Vert\cdot\right\Vert$ denotes the norm of the matrix acting on $\mathbb{C}^n$. It is worth observing that if $W$ is an $A_2$ matrix weight, then $W^{-1}$ is also an $A_2$ matrix weight.
 These boundedness results provided nontrivial matrix analogues of classical scalar results and spurred a wave of interest in harmonic analysis in the matrix-weighted and more general settings.  See, for example \cite{bw15b, bpw16, ct15, gol02, gol03, IKP,  k2, lt07, nptv02, nr15, p07, ps12, vt97b}.

Recall that in the scalar setting, $A_2$ weights $w$ are positive a.e., locally integrable functions satisfying
\[ [w]_{A_2} = \sup_{I} \left( \langle w \rangle_I \langle w^{-1} \rangle_I \right) < \infty.\]
where $\langle w \rangle_I=\frac{1}{|I|} \int_I w(x) dx$ and the supremum is over all intervals $I$. While many scalar results involving such $A_2$ weights have been extended to or disproved in the matrix weighted setting, there are still many open questions. Arguably the most famous is the Matrix $A_2$ conjecture.  Namely, a celebrated result by T. Hyt\"onen \cite{th12} proved the scalar $A_2$ conjecture, i.e.~that
\begin{equation} \label{eqn:linear} \| T \|_{L^2(w) \rightarrow L^2(w)} \lesssim [w]_{A_2}\end{equation}
for all Calder\'on-Zgymund operators $T$, where the implied constant does not depend on $w$. 
However in the matrix setting, despite substantial work, the sharp dependence of the operator norms on $[W]_{A_2}$ is unknown. Currently, the best known bound for any non-trivial operator is the bound $[W]^{\frac{3}{2}}_{A_2}$ for sparse operators, with proofs appearing in both \cite{bw15, IKP}.  

In this paper, we study particular examples of $A_2$ matrix weights. These examples are motivated by the scalar power weights, which are weights of the form $w(x) = a |x|^{\gamma}$, where $a, \gamma \in \mathbb{R}.$ It is well known that such a $w(x)$ is an $A_2$ weight if and only if $a$ is positive and   $-1 < \gamma <1$. See e.g.~Example $9.1.7$ in \cite{g09}. 
In the scalar setting, these are useful both as simple examples to check conjectures and because they provide the example showing that the linear bound \eqref{eqn:linear} is sharp \cite{pet07}. Our goal is to provide important examples of $A_2$ weights in the matrix setting, which may act as a testing ground for conjectures or be used to build interesting objects, such as systems of PDEs like the ones studied in \cite{im16}.

\subsection{Outline of Paper.} 
In this paper, we study two generalizations of power functions in the matrix setting. To do this, we require additional preliminary definitions and results related to the structure of matrices. These are presented in Section \ref{sec:facts}.

Then in Section \ref{sec:t1}, we examine matrix functions where each entry is a power of $|x|$. Specifically, we define \emph{Type 1 matrix power functions} to be matrix functions of the form
\begin{equation} \label{eqn:mpower} W(x) = \begin{bmatrix}
a_{11} |x|^{\gamma_{11}} & \dots & a_{1n} |x|^{\gamma_{1n}} \\
\vdots &\ddots & \vdots\\
a_{n1} |x|^{\gamma_{n1}}  & \dots & a_{nn} |x|^{\gamma_{nn}}
\end{bmatrix},\end{equation}
where each $a_{ij} \in \mathbb{C}$ and each $\gamma_{ij} \in \mathbb{R}.$ Given a Type $1$ matrix power function, we define its \emph{coefficient matrix} $A$ to be the matrix
\[ A \equiv \begin{bmatrix}
a_{11} & \dots & a_{1n}  \\
\vdots &\ddots & \vdots\\
a_{n1}  & \dots & a_{nn} 
\end{bmatrix}.\]
In Subsection \ref{ssec:pd}, we characterize when Type $1$ matrix functions $W(x)$ are positive definite for a.e.~$x$. In Subsection \ref{ssec:A2}, we establish our main result Theorem \ref{thm:A2}, which characterizes when Type $1$ matrix power functions are $A_2$ matrix weights. We use this result to generate several nontrivial $A_2$ matrix weights in Example \ref{ex:A2} and in Remark \ref{rem:multi}, point out that our arguments and results generalize to two kinds of Type $1$ matrix power functions in several variables.

In Section \ref{sec:t2}, we consider matrix functions whose eigenvalues are powers of $|x|$. Specifically, we define \emph{Type 2 matrix power functions} to be  matrix functions of the form
\begin{equation} \label{eqn:power2}
W(x) = U(x)\begin{bmatrix}
    \alpha_1|x|^{\gamma_1}  & \dots  & 0 \\
    \vdots &  \ddots & \vdots \\
    0  & \dots  & \alpha_n|x|^{\gamma_n}
\end{bmatrix}
U^\ast(x),\end{equation}
where each $\alpha_i \in \mathbb{C}$, each $\gamma_{i} \in \mathbb{R},$  the interior matrix is diagonal, and $U(x)$ is unitary for a.e.~$x$. In Subsection \ref{ssec:local}, we provide necessary and sufficient conditions for these matrix functions to be locally integrable. It is clear that they are positive definite a.e.~precisely when each $\alpha_i$ is positive. However, any necessary and sufficient conditions for Type $2$ matrix power functions to be $A_2$ matrix weights will likely be complicated and depend heavily on the structure of the unitary $U(x)$. Indeed, in Example \ref{ex:notA2}, we show that the supposedly-simple matrix weight
\[W(x) = \begin{bmatrix}
     \cos x & -\sin  x \\
     \sin x & \cos x \\
\end{bmatrix} \begin{bmatrix}
     |x|^{\gamma_1} & 0 \\
    0 & |x|^{\gamma_2} \\
\end{bmatrix}\begin{bmatrix}
     \cos x & \sin x \\
     -\sin x & \cos x \\
\end{bmatrix}\]
is almost never an $A_2$ matrix weight. A necessary requirement is that $\gamma_1=\gamma_2$, which turns $W(x)$ into the trivial weight  $|x|^{\gamma_1} I_{2\times 2}.$ In \cite{b81, b85}, Bloom studied a variant of this example for the case where $\gamma_2 = -\gamma_1$ and concluded that this weight is a ``good weight.'' As our result shows that even these supposedly-nice weights are almost never $A_2$ matrix weights, we conjecture that Type $2$ power functions are rarely $A_2$ matrix weights.  Lastly, in Remark \ref{rem:extension2}, we point out that this (non)example generalizes to higher dimensions. 

\section*{Acknowledgements} The first author would like to thank Brett Wick for many valuable discussions about $A_2$ matrix weights and Joshua Isralowitz for insightful conversations related to $A_2$ matrix weights in several variables and their applications.

\section{Preliminary Matrix Theory} \label{sec:facts}

In order to study matrix power functions, we require standard facts about matrices, including when a matrix is positive definite (i.e.~self-adjoint with positive eigenvalues) and methods to compute the inverse of a matrix. The needed definitions and facts are presented below and  can be found in the monographs \cite{bh, hj}.  

Throughout this section, let $A$ be an $n \times n$ matrix with entries $a_{ij} \in \mathbb{C}.$ The needed characterization of positive definiteness requires submatrices. A \emph{submatrix} of $A$ is a matrix composed of the entries of $A$ that lie in a subcollection of its rows and columns. Specifically, let $\alpha \subseteq \{ 1, \dots, n\}$ and $\beta \subseteq \{1, \dots, n\}$ be index sets. The notation $A [\alpha, \beta]$ indicates the submatrix of $A$ composed of the entries of $A$ from the rows of $A$ indexed by $\alpha$ and the columns of $A$ indexed by $\beta.$ A submatrix $A [\alpha, \beta]$  is a \emph{principal submatrix} if $\alpha =\beta$.  A principle submatrix $A[\alpha, \alpha]$ is a \emph{leading principal submatrix} if there is a $k \in \{1, \dots, n\}$ such that $\alpha=\{1, \dots, k\}.$ Finally, the shorthand notation $A_{ij} = A[ \{i\}^c, \{j\}^c]$ represents the submatrix of $A$ obtained by removing the i$^{th}$ row and j$^{th}$ column from $A$. 

The determinant of a (principle, leading principle) submatrix of $A$ is called a (\emph{principle}, \emph{leading principle}) \emph{minor} of $A$. One can determine if a self-adjoint matrix is positive definite using its principle and leading principle minors, via a characterization called \emph{Sylvester's Criterion}:

\begin{theorem} \label{thm:pos} Let $A$ be an $n \times n$ self-adjoint matrix. Then
\begin{itemize}
\item[a.] $A$ is positive definite if and only if every principal minor of $A$ is positive. 
\item[b.] $A$ is positive definite if and only if every leading principal minor of $A$ is positive.
\end{itemize}
\end{theorem}

To study $A_2$ matrix weights, we also need information about matrix inverses. This requires formulas for determinants, which require permutations. 
A permutation of length $n$ is a one-to-one function $\sigma: \{1,\dots, n\} \rightarrow \{1, \dots, n\},$ and the sign of $\sigma,$ denoted $\text{sgn}(\sigma),$ is $+1$ or $-1$ depending on whether the minimum number of transpositions needed to turn $\{1,\dots, n\}$ into $\sigma$ is even or odd.  Lastly, $S_n$ denotes the set of all permutations of length $n$. Then the determinant of $A$ can be computed using the formula
\begin{equation} \label{eqn:detA} \det A = \sum_{\sigma \in S_n} \left(  \text{sgn}(\sigma) \prod_{k=1}^n  a_{k \sigma(k)} \right).\end{equation}
If $A$ is invertible, we  can compute $A^{-1}$ using minors and determinants as follows
\begin{equation} 
\label{eqn:inverse}
A^{-1} = \frac{1}{\det A}\left[
\begin{array}{ccc} 
C_{11} & \dots & C_{n1} \\
\vdots &  \ddots & \vdots \\
C_{1n}  & \dots & C_{nn} \\
\end{array} \right], 
\end{equation}
where $C_{ij}$ is the $ij^{th}$  cofactor of $A$, meaning $C_{ij} = (-1)^{i+j}  \det A_{ij}$, where $A_{ij}$ was defined earlier. Note that $A^{-1}$ is $\frac{1}{\det A}$ times the transpose of the matrix of cofactors of $A$.

\section{Type $1$ Matrix Power Functions} \label{sec:t1}

In this section, we study Type $1$ matrix power functions, i.e.~matrix functions of form \eqref{eqn:mpower} with coefficient matrix $A$. 

\subsection{Positive Matrix Power Functions.} \label{ssec:pd} We will first characterize when a Type $1$ matrix power function $W(x)$ is positive definite a.e.
\begin{theorem} \label{thm:pd} An $n\times n$ Type $1$ matrix power function $W(x)$
is positive definite a.e.~if and only if its coefficient matrix $A$ is positive definite and 
the powers $\gamma_{ij}$ satisfy
\begin{equation}\label{eqn:gamma} \gamma_{ij} = \frac{\gamma_{ii} + \gamma_{jj}}{2} \ \text{ for } 1 \le i,j \le n.\end{equation}
\end{theorem}

This characterization requires the following lemma about determinant formulas of matrix power functions, which is proved in Subsection \ref{ssec:lemma}.

\begin{lemma} \label{lem:det} Let $W$ be an $n\times n$ Type $1$ matrix power function with coefficient matrix $A$ 
and powers $\gamma_{ij}$ satisfying $\gamma_{ij} = \frac{\gamma_{ii} + \gamma_{jj}}{2} \ \text{ for } 1 \le i,j \le n$. 
Let $W(x)_{ij}$ be the matrix obtained by removing the $i^{th}$ row and $j^{th}$ column from $W(x)$ for each $x$. Then 
\[ \det W(x) = |x|^{\sum_{k=1}^n \gamma_{kk}} \det A \ \ \text{ and } \ \ \det W(x)_{ij} = |x|^{ - \gamma_{ij} +\sum_{k=1}^n \gamma_{kk}}  \det A_{ij}.\] 
\end{lemma}

Here is the proof of Theorem \ref{thm:pd}:

\begin{proof} 
We will prove this by induction on the size of $W$.  For the base case, let $n=1.$ Then $W(x) = [a_{11}|x|^{\gamma_{11}}]$. 
 First assume $W(x)$ is positive definite a.e. Then since $|x|^{\gamma_{11}}> 0$ for $x\ne 0$, we can conclude $a_{11} > 0$. Thus, $A \equiv [a_{11}]$ is positive definite. In addition, condition \eqref{eqn:gamma} is trivial since $\frac{\gamma_{11} + \gamma_{11}}{2} = \gamma_{11}.$ 
Now assume $A = [a_{11}]$ is positive definite.  As $|x|^{\gamma_{11}}> 0$ for $x\ne 0$, this immediately implies $W(x) = [a_{11}|x|^{\gamma_{11}}]$ is positive definite. 
 
For the inductive step, assume the theorem holds for every $(n-1) \times (n-1)$ Type $1$ matrix power function, and let $W$ be an $n\times n$ Type $1$ matrix power function with coefficient matrix $A$ and powers $\gamma_{ij}.$ \\

\noindent ($\Rightarrow$) For the forward direction, assume $W(x)$ is positive definite a.e. We will show $A$ is positive definite and the $\gamma_{ij}$ powers satisfy condition \eqref{eqn:gamma}. First, since $W(x)$ is positive definite $a.e.$, it is self-adjoint a.e., i.e.
\[a_{ij}|x|^{\gamma_{ij}} = \overline{a_{ji}|x|^{\gamma_{ji}}} = \bar{a}_{ji} |x|^{\gamma_{ji}}, \qquad \text{ for } 1 \le i,j \le n.\]
This implies each $\gamma_{ij}=\gamma_{ji}$ and $a_{ij} = \bar{a}_{ji}$, so $A$ is self-adjoint. Furthermore, by Theorem \ref{thm:pos}, all of the leading principal minors of $W(x)$ are positive a.e. Specifically, for each $k=1, \dots, n-1$, 
\[ \det \left[
\begin{array}{ccc} 
a_{11} |x|^{\gamma_{11}} & \dots & a_{1k} |x|^{\gamma_{1k}} \\
\vdots & \ddots  & \vdots\\
a_{k1} |x|^{\gamma_{k1}}  & \dots & a_{kk} |x|^{\gamma_{kk}}
\end{array} \right]>0. \]
Thus, by Theorem \ref{thm:pos}, we can conclude that the $(n-1) \times (n-1)$  Type $1$ matrix power function
\[ W(x)_{nn} \equiv \left[
\begin{array}{ccc} 
a_{11} |x|^{\gamma_{11}} & \dots & a_{1(n-1)} |x|^{\gamma_{1(n-1)}} \\
\vdots & \ddots & \vdots\\
a_{(n-1)1} |x|^{\gamma_{(n-1)1}}  & \dots & a_{(n-1)(n-1)} |x|^{\gamma_{(n-1)(n-1)}}
\end{array} \right], \]
obtained by removing the $n^{th}$ row and $n^{th}$ column from $W(x),$ is also positive definite a.e. 
 Now, we establish condition \eqref{eqn:gamma}. First, applying the inductive hypothesis to $W(x)_{nn}$  gives
\[\gamma_{ij} = \frac{\gamma_{ii}+\gamma_{jj}}{2}, \qquad \text{ for } 
1 \leq i,j \leq n-1.\]
To conclude  \eqref{eqn:gamma}, we only need 
\[\gamma_{ni} =\gamma_{in} =  \frac{\gamma_{ii} + \gamma_{nn}}{2}, \qquad \text{ for } 1 \leq i \leq n,\] 
where the first equality was already established. Also, it is clear that $\gamma_{nn} = \frac{\gamma_{nn} + \gamma_{nn}}{2}$. Thus, we must show that $\gamma_{in} =  \frac{\gamma_{ii} + \gamma_{nn}}{2}$ for $1 \leq i \leq n-1$. Fix such an $i$. Then the following determinant
\[ \det \left[ \begin{array}{cc}
a_{ii}|x|^{\gamma_{ii}} & a_{in}|x|^{\gamma_{in}} \\
a_{ni}|x|^{\gamma_{ni}} & a_{nn}|x|^{\gamma_{nn}}
\end{array} \right] \]
is a principal minor of $W(x)$ and so by Theorem \ref{thm:pos}, is positive a.e. Computing this determinant and using the fact that $A$ is self-adjoint gives
\begin{equation} \label{eqn:posdet} a_{ii}a_{nn}|x|^{\gamma_{ii} + \gamma_{nn}} - |a_{in}|^2|x|^{2\gamma_{in}} > 0.\end{equation}
Here we can assume $a_{in} \ne 0$, because otherwise we could trivially choose $\gamma_{in} = \frac{\gamma_{ii}+\gamma_{nn}}{2}.$
Now by looking at both $x$ values near zero and $x$  values arbitrarily large, one can see that \eqref{eqn:posdet} holds a.e.~if and only if
\[ \gamma_{ii} + \gamma_{nn} = 2 \gamma_{in},\]
which is the desired equality. 
 
Now we show $A$ is positive definite. By the inductive hypothesis, the  coefficient matrix of $W(x)_{nn}$
\[ A_{nn} \equiv \left[\begin{array}{ccc} 
a_{11}  & \dots  & a_{1(n-1)}  \\
\vdots  & \ddots  & \vdots  \\
a_{(n-1)1}  & \dots  & a_{(n-1)(n-1)} 
\end{array} \right] \]
is positive definite a.e. 
Then Theorem \ref{thm:pos} implies that the leading principal minors of $A_{nn}$ are positive. Namely, 
\[ \det \left[\begin{array}{ccc} 
a_{11}  & \dots  & a_{1k}  \\
\vdots  & \ddots  & \vdots  \\
a_{k1}  & \dots  & a_{kk} 
\end{array} \right] >0,\]
for $k=1, \dots, n-1.$ Now, we show that $\det A > 0$. As $W(x)$ is positive definite a.e., $\det W(x) >0$ a.e. As Lemma \ref{lem:det} gives $\det W(x) = |x|^{\sum_{k=1}^n \gamma_{kk}} \det A$, we can conclude $\det A >0.$ This shows that all of the leading principal minors of $A$ are positive, so $A$ is positive definite.\\


\noindent ($\Leftarrow$) Now assume $A$ is positive definite and condition \eqref{eqn:gamma} holds. As $A$ is self-adjoint, $a_{ij} = \bar{a}_{ji}$ for $1 \leq i,j \leq n$. As condition \eqref{eqn:gamma} implies that each $\gamma_{ij} = \gamma_{ji}$, $W(x)$ is also self-adjoint. 
Using Theorem \ref{thm:pos} and the positive definiteness of $A$, one can easily show that the $(n-1) \times (n-1)$ matrix 
\[A_{nn} \equiv \left[\begin{array}{ccc} 
a_{11}  & \dots  & a_{1(n-1)}  \\
\vdots  &  \ddots & \vdots  \\
a_{(n-1)1}  & \dots  & a_{(n-1)(n-1)} 
\end{array} \right] \]
is positive definite. 
Then by the inductive hypothesis, the $(n-1) \times (n-1)$ Type $1$ matrix power function 
\[ W(x)_{nn} \equiv \left[
\begin{array}{ccc} 
a_{11} |x|^{\gamma_{11}} & \dots & a_{1(n-1)} |x|^{\gamma_{1(n-1)}} \\
\vdots & \ddots & \vdots\\
a_{n1} |x|^{\gamma_{(n-1)1}}  & \dots & a_{(n-1)(n-1)} |x|^{\gamma_{(n-1)(n-1)}}
\end{array} \right] \]
is positive definite a.e. Again by Theorem \ref{thm:pos}, this implies that all of the leading principal minors of $W(x)$ of the form
\[ \det \left[
\begin{array}{ccc} 
a_{11} |x|^{\gamma_{11}} & \dots & a_{1k} |x|^{\gamma_{1k}} \\
\vdots & \ddots & \vdots\\
a_{k1} |x|^{\gamma_{k1}}  & \dots & a_{kk} |x|^{\gamma_{kk}}
\end{array} \right] >0, \]
a.e.~for $k=1, \dots, n-1.$ Furthermore, since $\det A >0$ and by Lemma \ref{lem:det}, $\det W(x) = |x|^{\sum_{k=1}^n \gamma_{kk}} \det A$, we also have 
$\det W(x)> 0$ a.e. Thus, all leading principal minors of $W(x)$ are positive a.e., which implies $W(x)$ is positive definite a.e.~and completes the proof. \end{proof}
%

\subsection{Matrix $A_2$ Power Weights.} \label{ssec:A2} Now, we characterize when Type $1$ matrix power functions are actually $A_2$ matrix weights. 

\begin{theorem} \label{thm:A2}
Let  $W$ be a Type 1 $n\times n$ matrix power function with coefficient matrix $A$ and powers $\gamma_{ij}.$ Then $W$ is an $A_2$ matrix weight if and only if $A$ is positive definite and the $\gamma_{ij}$ satisfy
\[  \gamma_{ij} = \frac{\gamma_{ii} + \gamma_{jj}}{2}  \ \text{ and } \ -1 < \gamma_{ii} < 1  \qquad  \text{ for }  1 \leq i,j \leq n.\] 
\end{theorem}

To establish this, we require the following lemma, which is proved in Subsection \ref{ssec:lemma}.

\begin{lemma} \label{lem:inverse} Let  $W$ be a Type 1 $n\times n$ matrix power function with coefficient matrix $A$ and powers $\gamma_{ij}$ that is positive definite a.e. Then
\[ W(x)^{-1} = \frac{1}{\det A}\left[
\begin{array}{ccc} 
c_{11} |x|^{-\gamma_{11}} & \dots & c_{n1} |x|^{-\gamma_{1n}}\\
\vdots & \ddots & \vdots \\
c_{1n} |x|^{-\gamma_{n1}} & \dots & c_{nn}  |x|^{-\gamma_{nn}}\\
\end{array} \right],\]
where each $c_{ij} = (-1)^{i+j} \det A_{ij}.$
\end{lemma}

Now we prove Theorem \ref{thm:A2}:

\begin{proof}

($\Rightarrow$)
First, assume that $W$ is an $A_2$ matrix weight. Then $W(x)$ is positive definite a.e.~and so by Theorem \ref{thm:pd}, $A$ is positive definite  and 
\[ \gamma_{ij} = \frac{\gamma_{ii} + \gamma_{jj}}{2} \qquad  \text{ for } 1 \leq i,j \leq n.\]
To complete this direction of the proof, we just need to establish the bounds on the $\gamma_{ii}.$ 

Observe that the diagonal entries  $a_{ii}|x|^{\gamma_{ii}}$ of $W(x)$ are exactly the $1 \times 1$ principle minors of $W(x).$ Then, as $W(x)$ is positive definite a.e., Theorem \ref{thm:pos} implies that each $a_{ii}|x|^{\gamma_{ii}}$ is positive $a.e.$ Thus, each $a_{ii}$ is nonzero. Furthermore, by Lemma \ref{lem:inverse},  
\[ W(x)^{-1} = \frac{1}{\det A}\left[
\begin{array}{ccc} 
c_{11} |x|^{-\gamma_{11}} & \dots & c_{n1} |x|^{-\gamma_{1n}}\\
\vdots & \ddots & \vdots \\
c_{1n} |x|^{-\gamma_{n1}} & \dots & c_{nn}  |x|^{-\gamma_{nn}}\\
\end{array} \right],\]
where each $c_{ij} = (-1)^{i+j} \det A_{ij}.$ As $W(x)^{-1}$ is positive definite a.e., we can apply Theorem \ref{thm:pos} to conclude that each $c_{ii}$ is nonzero. 

Since $W$ and hence, $W^{-1}$ are $A_2$ matrix weights, their entries are locally integrable. In particular, the diagonal entries $a_{ii} |x|^{\gamma_{ii}}$  and $c_{ii} |x|^{-\gamma_{ii}}$ are integrable on intervals containing the origin. As each $a_{ii}$ and $c_{ii}$ are nonzero, this implies $ -1 < \gamma_{ii} <1,$ as desired.\\

\noindent ($\Leftarrow$)
Now, assume that $A$ is positive definite and the $\gamma_{ij}$ satisfy the given conditions. Then by Theorem \ref{thm:pd},  $W(x)$ is positive definite a.e.  Moreover, the conditions on the $\gamma_{ij}$ imply that  $-1 < \gamma_{ij} < 1$ for $1 \le i,j\le n,$  and so each entry of $W$ is locally integrable. Thus, $W$ is a matrix weight.

To show $W$ is an $A_2$ matrix weight, we need both $W$ and $W^{-1}.$ Using the given formula for $W$ and the one for $W^{-1}$ in Lemma \ref{lem:inverse}, we can compute their averages over intervals componentwise. Indeed for any interval $I$,
\[ \left \langle W \right \rangle_I =   \left[
\begin{array}{ccc} 
a_{11} \langle |x|^{\gamma_{11}} \rangle_I & \dots &  a_{1n} \langle|x|^{\gamma_{1n}} \rangle_I \\
\vdots & \ddots  & \vdots \\
a_{n1} \langle |x|^{\gamma_{n1}} \rangle_I & \dots &a_{nn} \langle |x|^{\gamma_{nn}} \rangle_I\\
\end{array} \right] \]
and similarly
\[ \left \langle W^{-1} \right \rangle_I  = \frac{1}{\det A}\left[
\begin{array}{ccc} 
c_{11} \langle |x|^{-\gamma_{11}} \rangle_I & \dots & c_{n1} \langle |x|^{-\gamma_{1n}}\rangle_I \\
\vdots & \ddots & \vdots \\
c_{1n} \langle |x|^{-\gamma_{n1}} \rangle_I & \dots & c_{nn} \langle  |x|^{-\gamma_{nn}} \rangle_I \\
\end{array} \right].
\]
To show that $W$ is an $A_2$ matrix weight, we simply need to verify that
\[ [W]_{A_2} \approx \sup_I \text{Tr} \left(\langle W \rangle_I \langle W^{-1} \rangle_I \right)< \infty.\]
To see this, fix an interval $I$ and observe that
\[ 
\text{Tr}(\langle W \rangle_I \langle W^{-1} \rangle_I) =  \sum_{i,j=1}^n \left( \langle W \rangle_I\right)_{ij} \left( \langle W^{-1} \rangle_I\right)_{ji} 
\leq \displaystyle  \frac{1}{|\det A|} \sum_{i,j = 1}^n |a_{ij}c_{ij}| \langle |x|^{\gamma_{ij}} \rangle_I \langle |x|^{-\gamma_{ij}} \rangle_I. 
\]
Notice that each $\gamma_{ij}$ satisfies $-1 < \gamma_{ij}<1,$ so $|x|^{\gamma_{ij}}$ is a scalar $A_2$ weight. Thus, we can conclude that
\[ \text{Tr}(\langle W \rangle_I \langle W^{-1} \rangle_I) \le  \sup_{i,j} \left( \frac{|a_{ij} \det A_{ij}|} {|\det A|}  \right) \displaystyle \sum_{i,j=1}^n   \big[|x|^{\gamma_{ij}} \big]_{A_2},\]
which is bounded independent of $I$. Hence, $W$ is an $A_2$ matrix weight. 
\end{proof}

Alternately, one could prove that these positive definite power functions are matrix $A_2$ weights using the sufficient conditions from Theorem $ 4.3$ in \cite{nr15} or by proving they are ``almost diagonal'' and using Proposition $4.2$ in \cite{b81}.  

Now, let us apply this theorem to generate nontrivial examples of $A_2$ matrix weights.

\begin{example} \label{ex:A2} To obtain an $n \times n$ Type $1$ matrix $A_2$ power weight, one need only choose a positive coefficient matrix $A$ and diagonal powers $-1< \gamma_{11}, \dots, \gamma_{nn}< 1$. For example, consider the positive definite matrix
\[ A = \begin{bmatrix*}[r] 5 & 3 \\ 3 & 2  \end{bmatrix*}
\ \ 
\text{ and powers } \ \  \gamma_{11} = \frac{1}{2}, \gamma_{22} = -\frac{2}{3}. \]
These generate the $A_2$ matrix weight
\[ W(x) =  \begin{bmatrix*}[c] 5 |x|^{1/2} & 3|x|^{-1/12} \\ 3|x|^{-1/12}  & 2 |x|^{-2/3} \end{bmatrix*}.\]
Similarly, we can use  the positive definite matrix
\[
B = \begin{bmatrix*}[r] 4 & 1 & 2 \\
1 & 2 & -1 \\
2 & -1 & 3 \end{bmatrix*} \ \ \text{ and powers } \ \ \gamma_{11} = \frac{3}{4}, \gamma_{22} = {-\frac{3}{4}}, \gamma_{33} = \frac{1}{2} \]
to generate the $A_2$ matrix weight
\[ W(x) = \begin{bmatrix*}[c] 4 |x|^{3/4} & 1 & 2|x|^{5/8} \\
1 & 2|x|^{-3/4} & -|x|^{-1/8} \\
2|x|^{5/8} & -|x|^{-1/8} & 3|x|^{1/2} \end{bmatrix*}.\]
\end{example}

Many of our arguments generalize almost immediately to two kinds of Type $1$ matrix power weights in several variables, as detailed below.

\begin{remark} \label{rem:multi} Here, we consider the $d$-variable setting, where $x = (x_1, \dots, x_d) \in \mathbb{R}^d.$ First, a matrix function $W$ in $d$ variables is an $A_2$ matrix weight if its entries are locally integrable and if it is positive definite a.e.~and satisfies
\[\big [ W \big ] _{A_2} \equiv \sup_{Q} \left \| \langle W \rangle_Q^{\frac{1}{2}}
 \langle W^{-1} \rangle_Q^{\frac{1}{2}} \right \|^2 < \infty,
\]
where the supremum is over cubes $Q = I_1 \times \cdots \times I_d$, with each pair of intervals $I_i$, $I_j$ satisfying $|I_i| = |I_j|.$ Now, we restrict to two variables to simplify notation. In two variables, define the Type $1.a$ matrix power functions to be matrix functions of the form 
\[ W(x) = \begin{bmatrix}
a_{11} |x_1|^{\gamma_{11}} |x_2|^{\beta_{11}}   & \dots & a_{1n} |x_1|^{\gamma_{1n}} |x_2|^{\beta_{1n}} \\
\vdots &\ddots & \vdots\\
a_{n1} |x_1|^{\gamma_{n1}} |x_2|^{\beta_{n1}}  & \dots & a_{nn} |x_1|^{\gamma_{nn}} |x_2|^{\beta_{nn}}
\end{bmatrix},
\]
with coefficient matrix $A$ and powers $\gamma_{ij}$. Without too much effort, one can generalize Lemmas \ref{lem:det} and \ref{lem:inverse} as well as Theorems \ref{thm:pd} and \ref{thm:A2} to this two variable setting.  
In this setting, the main result states:~a Type 1.a $n\times n$ matrix power function $W$ with coefficient matrix $A$ and powers $\gamma_{ij}$ is an $A_2$ matrix weight if and only if $A$ is positive definite and the $\gamma_{ij}$ and $\beta_{ij}$ satisfy
\[  \gamma_{ij} = \frac{\gamma_{ii} + \gamma_{jj}}{2}, \beta_{ij} = \frac{\beta_{ii} + \beta_{jj}}{2}   \ \text{ and } \ -1 < \gamma_{ii}, \beta_{ii}  < 1  \qquad  \text{ for }  1 \leq i,j \leq n.\] 
The obvious generalization holds in $d$ variables.

Similarly, define $d$ variable Type $1.b$ matrix power functions to be matrix functions of the form
\[ W(x) = \begin{bmatrix}
a_{11} \|x\|^{\gamma_{11}}  & \dots & a_{1n} \|x\|^{\gamma_{1n}} \\
\vdots &\ddots & \vdots\\
a_{n1} \|x\|^{\gamma_{n1}}  & \dots & a_{nn} \|x\|^{\gamma_{nn}} 
\end{bmatrix},
\]
with coefficient matrix $A$, where $ \| x \| = \sqrt { x_1^2 + \cdots + x_d^2}.$ 
 One can also generalize Lemmas \ref{lem:det} and \ref{lem:inverse} as well as Theorems \ref{thm:pd} and \ref{thm:A2} to these types of matrix functions. Here the generalization of Theorem \ref{thm:A2} is slightly different, as $\| x \|^{\gamma}$ is locally integrable if and only if $\gamma > -d$ and $\| x \|^{\gamma}$ is an $A_2$ weight if and only if $-d < \gamma < d.$ For this fact, see Example $9.1.7$ in \cite{g09}. Then in this setting, the main result states:~a Type 1.b $n\times n$ matrix power function  $W$ with coefficient matrix $A$ and powers $\gamma_{ij}$ is an $A_2$ matrix weight if and only if $A$ is positive definite and the $\gamma_{ij}$ satisfy
\[  \gamma_{ij} = \frac{\gamma_{ii} + \gamma_{jj}}{2}    \ \text{ and } \ -d < \gamma_{ii} < d  \qquad  \text{ for }  1 \leq i,j \leq n.\]

\end{remark}

\subsection{Proofs of Lemmas \ref{lem:det} and \ref{lem:inverse}} \label{ssec:lemma} 
 In this subsection, we include the proofs of the lemmas used earlier.  Here is the proof of Lemma \ref{lem:det}:

\begin{proof} First we show $\det W(x) = |x|^{\sum_{k=1}^n \gamma_{kk}} \det A$. Let $(W(x))_{ij}$ denote the $ij^{th}$ entry of $W(x).$ Then by \eqref{eqn:detA}, we have
\begin{equation} \label{eqn:det2} \det W(x) = \sum_{\sigma \in S_n} \left(\text{sgn}(\sigma) \prod_{k=1}^n (W(x))_{k\sigma(k)}\right).\end{equation}
Fix $\sigma \in S_n$ and consider $\prod_{k=1}^n (W(x))_{k\sigma(k)}$. As $k = 1,...,n$ and $\sigma$ is injective, the product has one entry from each row and column of $W(x)$. Then
 
\begin{align*}\prod_{k=1}^n (W(x))_{k\sigma(k)} &= \prod_{k=1}^n a_{k\sigma(k)}|x|^{\gamma_{k\sigma(k)}}\\
 &= \prod_{k=1}^n a_{k\sigma(k)}|x|^{(\gamma_{kk} + \gamma_{\sigma(k)\sigma(k)})/2}\\ 
&= \prod_{k=1}^n a_{k\sigma(k)} \prod_{k=1}^n |x|^{\gamma_{kk}/2} \prod_{\sigma(k)=1}^n |x|^{\gamma_{\sigma(k)\sigma(k)}/2}\\
&= \prod_{k=1}^n a_{k\sigma(k)} |x|^{\sum_{k=1}^n \gamma_{kk}}.
\end{align*} 
By substituting this into \eqref{eqn:det2} and factoring  $|x|^{\sum_{k=1}^n \gamma_{kk}} $ out of the sum, we obtain
\begin{align*}\det W(x) = \sum_{\sigma \in S_n} \left(\text{sgn}(\sigma)  |x|^{\sum_{k=1}^n \gamma_{kk}} \prod_{k=1}^n a_{k\sigma(k)}\right) =  |x|^{\sum_{k=1}^n \gamma_{kk}}\det A,\end{align*}
which  gives the first equation.
Now we fix $i,j$ and show 
\[ \det W(x)_{ij} =  |x|^{ - \gamma_{ij} +\sum_{k=1}^n \gamma_{kk}}\det A_{ij}.\]
 First by \eqref{eqn:detA}, we have
\begin{equation} \label{eqn:det3} \det W(x)_{ij} = \sum_{\sigma \in S_{n-1}} \left(\text{sgn}(\sigma) \prod_{k=1}^{n-1}(W(x)_{ij})_{k\sigma(k)} \right).\end{equation} 
Fix $\sigma \in S_{n-1}$ and consider $\prod_{k=1}^{n-1} (W(x)_{ij})_{k\sigma(k)}$. 
As $k = 1,...,n-1$ and $\sigma$ is injective, the product has one entry from each row and column of $W(x)_{ij}$. Define the indices
\[\beta(k) \equiv \left \{ \begin{array}{ll} k & \text{ if } k < i, \\ 
k+1 & \text{ if } k \geq i, \end{array} \right.
\ \ \text{and}  \ \ \phi(k) \equiv \left \{ \begin{array}{ll} \sigma(k) & \text{ if } \sigma(k) < j, \\ 
\sigma(k)+1 & \text{ if } \sigma(k) \geq j. \end{array} \right.\] 
Then using arguments similar to those used to obtain the previous determinant formula, 
\[ \prod_{k=1}^{n-1} (W(x)_{ij})_{k\sigma(k)}   =  \prod_{k=1}^{n-1} (A_{ij})_{k\sigma(k)}|x|^{\gamma_{\beta(k)\phi(k)}}
=  |x|^{ - \gamma_{ij} +\sum_{k=1}^n \gamma_{kk}} \prod_{k=1}^{n-1} (A_{ij})_{k\sigma(k)}.\]
Substituting that into \eqref{eqn:det3} gives the desired formula.\end{proof}

Here is the proof of Lemma \ref{lem:inverse}: 

\begin{proof} Using the inverse formula \eqref{eqn:inverse}, we have
\[ W(x)^{-1} \equiv \frac{1}{\det W(x)}\left[
\begin{array}{ccc} 
C_{11}(x) & \dots & C_{n1}(x) \\
\vdots &  \ddots & \vdots \\
C_{1n}(x)  & \dots & C_{nn}(x) \\
\end{array} \right], 
\]
where $C_{ij}(x) = (-1)^{i+j}  \det W(x)_{ij}$ is the $ij^{th}$  cofactor of $W(x)$. As $W(x)$ is positive definite a.e., we can use Theorem \ref{thm:pd} to conclude that 
$W$ satisfies the conditions of  Lemma \ref{lem:det}. Then by Lemma \ref{lem:det}, we know that each
\[  C_{ij}(x) = (-1)^{i+j} |x|^{ -\gamma_{ij} + \sum_{k=1}^n \gamma_{kk}} \det A_{ij} \ \
\text{ and}   \ \ \det W(x) =  |x|^{\sum_{k=1}^n \gamma_{kk}} \det A.\]
To simplify notation, define  $c_{ij} = (-1)^{i+j} \det A_{ij}$ for each $ 1 \le i,j\le n.$  Then combining formulas gives
\[ 
\begin{aligned}
W(x)^{-1} &\equiv \frac{{|x|^{-\sum_{k = 1}^n \gamma_{kk}}}}{ \det A}\left[
\begin{array}{ccc} 
c_{11} |x|^{-\gamma_{11} + \sum_{k=1}^n \gamma_{kk}}  & \dots & c_{n1} |x|^{-\gamma_{n1} + \sum_{k=1}^n \gamma_{kk}} \\
\vdots & \ddots & \vdots \\
 c_{1n} |x|^{-\gamma_{1n} + \sum_{k=1}^n \gamma_{kk}} & \dots & c_{nn} |x|^{-\gamma_{nn} + \sum_{k=1}^n \gamma_{kk}} \\
\end{array} \right]  \\
&= \frac{1}{\det A}\left[
\begin{array}{ccc} 
c_{11} |x|^{-\gamma_{11}} & \dots & c_{n1} |x|^{-\gamma_{1n}}\\
\vdots & \ddots & \vdots \\
c_{1n} |x|^{-\gamma_{n1}} & \dots & c_{nn}  |x|^{-\gamma_{nn}}\\
\end{array} \right],
\end{aligned}
\]
as desired. \end{proof}

\section{Type $2$ Matrix Power Functions} \label{sec:t2}

In this section, we consider Type $2$ matrix power functions, which have form \eqref{eqn:power2}. Our original goal was to characterize when such $W$ are $A_2$ matrix weights, but this appears to depend very closely on the structure of the unitary $U(x)$. Nevertheless, we have determined some necessary conditions for these matrix functions to be $A_2$ matrix weights.

\subsection{Locally Integrable Matrix Power Functions.} \label{ssec:local}

If $W$ is an $A_2$ matrix weight, then $W(x)$ and $W(x)^{-1}$ must be positive definite a.e.~and have locally integrable entries.  For Type $2$ matrix power functions, it is straightforward to characterize when this occurs.

\begin{theorem} \label{thm:local} Let $W$ be an  $n \times n$ Type $2$ matrix power function with unitary $U(x)$ and associated eigenvalues $\alpha_i |x|^{\gamma_i}$.  Then $W(x)$ is positive definite a.e.~and has locally integrable entries if and only if $\gamma_i > -1$ and $\alpha_i >0$ for  $ 1 \le i \le n$. \end{theorem}

\begin{proof} Let $u_{ij}(x)$ denote the $ij^{th}$ entry of $U(x)$ and let $w_{ij}(x)$ denote the $ij^{th}$ entry of $W(x)$. Then multiplying out the matrices in \eqref{eqn:power2} gives 
\[w_{ij}(x) = \sum_{k=1}^n \alpha_k|x|^{\gamma_k} u_{ik}(x)\overline{u_{jk}(x)}.\]
\noindent ($\Rightarrow$) Suppose $W(x)$ is positive definite a.e.~with locally integrable entries. Because the eigenvalues of $W(x)$ are $\alpha_i |x|^{\gamma_i}$, it follows that  each $\alpha_i$ is positive. By assumption, each $w_{ii}(x)$ is locally integrable, so the following sum is locally integrable
 \[\begin{aligned}\sum_{i=1}^n w_{ii}(x) &= \sum_{i=1}^n\sum_{k=1}^n \alpha_k|x|^{\gamma_k} u_{ik}(x)\overline{u_{ik}(x)}\\
&= \sum_{k=1}^n\alpha_k|x|^{\gamma_k}   \left( \sum_{i=1}^n |u_{ik}(x)|^2 \right) \\
&= \sum_{k=1}^n \alpha_k|x|^{\gamma_k}.\end{aligned}\]
Local integrability implies that
\begin{equation} \label{eqn:localint} \int_I \left(\sum_{k=1}^n \alpha_k|x|^{\gamma_k} \right) dx = \sum_{k=1}^n \alpha_k\int_I|x|^{\gamma_k}dx < \infty,\end{equation} for all finite intervals $I$.  Since each  $\alpha_k$ is positive, we can conclude that each $\gamma_k> -1$. Otherwise,  inequality \eqref{eqn:localint} would fail for intervals containing zero.\\

\noindent ($\Leftarrow$) Now, assume that each $ \gamma_i>-1$ and each $\alpha_i >0$. As $W(x)$ is self-adjoint and the eigenvalues of $W(x)$ are $\alpha_i |x|^{\gamma_i}$, it follows that $W(x)$ is positive definite a.e. To see that each $w_{ij}(x)$ is locally integrable, observe that
\[  |w_{ij}(x)| = \left |\sum_{k=1}^n \alpha_k|x|^{\gamma_k} u_{ik}(x)\overline{u_{jk}(x)} \right| \le \ \sum_{k=1}^n\alpha_k |x|^{\gamma_k}, \]
where we used the fact that since $U(x)$ is unitary, each $|u_{ik}(x) |\le 1.$ As each $ \gamma_i>-1$ implies each function in the final sum is locally integrable, each $w_{ij}(x)$ is locally integrable as well.
  \end{proof}

This result has implications for characterizing when Type $2$ matrix power functions are also $A_2$ matrix weights.

\begin{remark} \label{rem:a2} If $W$ is both a Type $2$ matrix power function and an $A_2$ matrix weight, then both $W$ and $W^{-1}$ are positive definite a.e.~with locally integrable entries. Then Theorem \ref{thm:local} implies that the coefficients $\alpha_i$ are positive and the powers $\gamma_i$ satisfy $-1< \gamma_i < 1.$ 
\end{remark}

It is worth noting that results similar to Theorem \ref{thm:local}, with more restrictive unitary matrices and more general eigenvalues, are discussed by Bloom in \cite{b81}. 

Now we  consider a specific example that shows that Type $2$ matrix power functions can satisfy the necessary conditions from Remark \ref{rem:a2} but fail to be $A_2$ matrix weights. In this example, the unitary $U(x)$ is the standard two dimensional rotation matrix.  In \cite{b81, b85}, Bloom studied a variant of this example when $\gamma_2 = -\gamma_1$ and determined that, by certain criteria, the weight is well behaved.  In contrast, we show that this matrix weight almost always fails to be an $A_2$ matrix weight.

\begin{example} \label{ex:notA2}
Define the matrix weight $$W(x) = \begin{bmatrix}
     \cos x & -\sin  x \\
     \sin x & \cos x \\
\end{bmatrix} \begin{bmatrix}
     |x|^{\gamma_1} & 0 \\
    0 & |x|^{\gamma_2} \\
\end{bmatrix}\begin{bmatrix}
     \cos x & \sin x \\
     -\sin x & \cos x \\
\end{bmatrix}.$$
Then $W$ is an $A_2$ matrix weight if and only if $\gamma_1 = \gamma_2$ and $-1 < \gamma_1 < 1$.

\end{example}

\begin{proof}

($\Rightarrow$) We prove this by contrapositive and assume  $\gamma_1 \ne \gamma_2$ or $\gamma_1$ does not satisfy the given inequality. Clearly, if $\gamma_1\ge1$ or $\gamma_1 \le -1$,  then Theorem \ref{thm:local} implies that $W$ does not have locally integrable entries and so, is not an $A_2$ matrix weight. Thus, we need only show that if $\gamma_1 \ne \gamma_2$, then $W$ is not an $A_2$ matrix weight.

By Lemma $3.6$ in \cite{tv97}, if $W$ is an $A_2$ matrix weight, then the diagonal entries of $W$ are scalar $A_2$ weights. So, to show that $W$ is not an $A_2$ matrix weight, we need only show that one of its diagonal entries, either $w_{11}(x)$ or $w_{22}(x)$, fails to be a scalar $A_2$ weight. Without loss of generality, assume $\gamma_2 > \gamma_1.$ In this situation, we claim  
$$w_{11}(x) = |x|^{\gamma_1}\cos^2x+|x|^{\gamma_2}\sin^2x$$
is not a scalar $A_2$ weight.  If instead $\gamma_1> \gamma_2$, one can apply identical arguments to show $w_{22}(x)$ is not a scalar $A_2$ weight. 

To show $w_{11}(x)$ is not an $A_2$ weight, define the interval $I_n = [2\pi n, 2\pi n + \pi]$ for each $n \in \mathbb{N}.$ We will show
 $$\lim_{n \rightarrow \infty} \left(\frac{1}{|I_n|}\int_{I_n} w_{11}(x) \ dx\right) \left(\frac{1}{|I_n|}\int_{I_n} \frac{1}{w_{11}(x)} \ dx \right)= \infty.$$
Choosing $n$ large enough, using the definition of $I_n$, and restricting to an interval where $\sin x \ge \frac{1}{\sqrt{2}}$ gives
\[\begin{aligned}\frac{1}{|I_n|}\int_{I_n} w_{11}(x) \ dx &= \frac{1}{\pi}\int_{2\pi n}^{2\pi n + \pi} w_{11}(x) \ dx \\
&\geq \frac{1}{\pi}\int_{2\pi n + \frac{\pi}{4}}^{2\pi n + \frac{3\pi}{4}} |x|^{\gamma_2}\sin^2x \ dx\\ 
&\geq \frac{1}{2\pi}\int_{2\pi n + \frac{\pi}{4}}^{2\pi n + \frac{3\pi}{4}} |x|^{\gamma_2} \ dx\\ 
&\gtrsim n^{\gamma_2},\end{aligned}\] 
where the implied constant does not depend on $n$.

Define $\epsilon_n = n^{(\gamma_1-\gamma_2)/2}$, where $n$ is large enough so $\epsilon_n < \pi$. For $n$ large enough, we always have $\epsilon_n<\pi$ because  $\gamma_1 - \gamma_2 < 0$. Furthermore, by the periodicity of $\sin x$ and the fact that $|\sin x | \le |x|$, we have
$$w_{11}(x) = |x|^{\gamma_1}\cos^2x+|x|^{\gamma_2}\sin^2x \leq |x|^{\gamma_1} + |x|^{\gamma_2}|x-2\pi n|^2.$$ Restricting $I_n$ to a smaller interval, we have 
\[\begin{aligned}\frac{1}{|I_n|}\int_{I_n} \frac{1}{w_{11}(x)} \ dx &\geq \frac{1}{\pi}\int_{2\pi n}^{2\pi n + \epsilon_n} \frac{dx}{|x|^{\gamma_1}\cos^2x+|x|^{\gamma_2}\sin^2x}\\
&\geq \frac{1}{\pi}\int_{2\pi n}^{2\pi n + \epsilon_n} \frac{dx}{|x|^{\gamma_1} + |x|^{\gamma_2}|x-2\pi n|^2}\\
&\approx \frac{1}{\pi}\int_{2\pi n}^{2\pi n + \epsilon_n} \frac{dx}{|2\pi n|^{\gamma_1} + |2\pi n |^{\gamma_2}\epsilon_n^2}\\
&\gtrsim \int_{2\pi n}^{2\pi n + \epsilon_n} \frac{dx}{n^{\gamma_1}+n^{\gamma_2}n^{\gamma_1-\gamma_2}}\\ 
&\approx n^{-\gamma_1}\epsilon_n\\ 
&= n^{\frac{-\gamma_1-\gamma_2}{2}},\end{aligned}\]
where the implied constant depends on $\gamma_1$ and $\gamma_2$ but  not on $n$. Thus for large $n$, 
$$\left(\frac{1}{|I_n|}\int_{I_n} w_{11}(x) \ dx\right) \left(\frac{1}{|I_n|}\int_{I_n} \frac{1}{w_{11}(x)} \ dx\right) \gtrsim (n^{\gamma_2})(n^{\frac{-\gamma_1-\gamma_2}{2}}) = n^{\frac{\gamma_2-\gamma_1}{2}}.$$
As $\gamma_2 - \gamma_1>0$, we can conclude $$\lim_{n\rightarrow\infty} n^{\frac{\gamma_2-\gamma_1}{2}} = \infty,$$ 
and so $w_{11}(x)$ is not a scalar $A_2$ weight. Thus, $W$ is not an $A_2$ matrix weight.\\

\noindent ($\Leftarrow$) For the other direction, assume $\gamma_1=\gamma_2$ and $-1 < \gamma_1 <1.$ Then $|x|^{\gamma_1}$ is a scalar $A_2$ weight and $ W(x) = |x|^{\gamma_1} I_{n\times n}$ is a matrix weight. Furthermore, for each interval $I$,
\[ \text{Tr }\left( \left \langle  W \right \rangle_I \left \langle W^{-1} \right \rangle_I \right) = \text{Tr}\left( \langle  |x|^{\gamma_1} \rangle_I \langle  |x|^{-\gamma_1} \rangle_I I_{n \times n} \right) \le n [|x|^{\gamma_1}]_{A_2}.\]
As this bound is independent of the interval $I$, $W$ is an $A_2$ matrix weight. \end{proof}

\begin{remark} \label{rem:extension2} \textnormal One can extend Example \ref{ex:notA2} to higher dimensions by observing that $W$ is built using a unitary matrix, which is a rotation matrix with angle $x$, and a diagonal matrix, whose entries are powers of $|x|$.  To extend to three dimensions, recall that the basic three dimensional rotation matrices are the following
\[ R_x(\theta) = 
\begin{bmatrix} 1 & 0 & 0 \\ 0 & \cos \theta & \sin \theta \\
0 &- \sin \theta & \cos \theta \end{bmatrix} \ \ 
R_y(\theta) = 
\begin{bmatrix} \cos \theta & 0 & -\sin \theta \\ 0 & 1 & 0 \\
\sin \theta & 0 & \cos \theta \end{bmatrix}
 \ \ 
R_z(\theta) = 
\begin{bmatrix} \cos \theta &  \sin \theta & 0 \\  -\sin \theta & \cos \theta & 0 \\
0 & 0 & 1 \end{bmatrix} 
 \]
where $R_{x}(\theta)$, $R_{y}(\theta)$ and $R_z(\theta)$ respectively rotate a vector around the $x$-, $y$-, or $z$-axis by an angle $\theta$. Every three dimensional rotation can be written as a product $R(\alpha, \beta, \gamma)\equiv R_x(\alpha)R_y(\beta) R_z(\gamma) $, where $\gamma, \beta, \alpha$ are the Euler angles of the rotation. More specifically, a general rotation matrix  in three dimensions can be represented as
\[R(\alpha, \beta, \gamma) = \begin{bmatrix} 
 \cos \beta \cos \gamma &   \cos \beta \sin \gamma & -\sin \beta & \\
 \sin \alpha  \sin \beta \cos \gamma -\cos \alpha  \sin \gamma  & \sin \alpha \sin \beta \sin \gamma + \cos \alpha  \cos \gamma&  \sin \alpha  \cos \beta   \\
 \cos \alpha \sin \beta \cos \gamma + \sin \alpha \sin \gamma  &\cos \alpha \sin \beta  \sin \gamma - \sin \alpha  \cos \gamma& \cos \alpha  \cos \beta
\end{bmatrix}.
\]
For details, see \cite{a11}, pp.~$59.$ To define a matrix function $W$, we need to specify a unitary and diagonal matrix. By setting $\alpha = \beta= \gamma = x$, we obtain the following three dimensional rotation matrix:
\[ U(x) \equiv\begin{bmatrix} 
\cos^2x &\cos x \sin x &  -\sin x \\
  \cos x \sin^2 x -  \cos x \sin x  & \cos^2x + \sin^3 x & \cos x \sin x   \\
\cos^2 \sin x + \sin^2 x &\cos x \sin^2x - \cos x\sin x  & \cos^2 x 
\end{bmatrix}.
\]
Mirroring Example \ref{ex:notA2}, define the Type $2$ matrix power function
\[ W(x) = U(x) \begin{bmatrix} |x|^{\gamma_1} & 0 & 0 \\ 0 & |x|^{\gamma_2} & 0 \\ 0 & 0 & |x|^{\gamma_3} \end{bmatrix} U^*(x). \]
Then, one can use arguments very similar to those in Example \ref{ex:notA2} to show that $W$ is an $A_2$ matrix weight if and only if $\gamma_2 = \gamma_2 = \gamma_3$ and $-1 < \gamma_1 <1.$

\end{remark}

\end{document}